\documentclass[11pt,twoside]{article}

\usepackage{amsbsy,amsfonts,amsmath,amssymb,eucal,mathrsfs,amsthm}
\usepackage[all]{xy}
\usepackage{pstricks}
\usepackage{pst-plot}
\usepackage{youngtab}

\newcommand{\GL}{\mathrm{GL}}
\newcommand{\SL}{\mathrm{SL}}
\newcommand{\Ker}{\mathrm{Ker}}

\newcommand{\Sp}{\mathrm{Sp}}
\newcommand{\SO}{\mathrm{SO}}
\newcommand{\Or}{\mathrm{O}}
\newcommand{\gl}{\mathfrak{gl}}

\newcommand{\so}{\mathfrak{so}}
\newcommand{\Y}{\mathfrak{Y}}

\theoremstyle{definition}
\newtheorem{definition}{Definition}[section]
\newenvironment{defi}{\begin{definition} \rm}{\end{definition}}

\newtheorem{notation}[definition]{Notation}
\newenvironment{nota}{\begin{notation} \rm}{\end{notation}}
\newtheorem{construction}[definition]{Elementary construction}
\newenvironment{cons}{\begin{construction} }{\end{construction}}
\newtheorem{example}[definition]{Example}

\newtheorem{examples}[definition]{Examples}

\newtheorem{nothing}[definition]{$\!\!$}

\theoremstyle{plain}
\newtheorem{prop}[definition]{Proposition}

\newtheorem{lemm}[definition]{Lemma}

\newtheorem{coro}[definition]{Corollary}
\newtheorem{theo}[definition]{Theorem}

\theoremstyle{remark}
\newtheorem{remark}[definition]{Remark}
\newenvironment{rema}{\begin{remark} \rm}{\end{remark}}
\newtheorem{remarks}[definition]{Remarks}

\newenvironment{proo}{{\flushleft \it Proof.}}{\hfill $\square$
  \vspace{2mm}}

\newenvironment{proo-norm}{{\flushleft \it Proof of the
    normality.}}{\hfill $\square$ \vspace{2mm}}
\newenvironment{proo-rat}{{\flushleft \it Proof of rationality of the
    resolution.}}{\hfill $\square$ \vspace{2mm}}

\theoremstyle{definition}
\newtheorem{definition*}{Definition}[section]
\newenvironment{defi*}{\begin{definition*} \rm}{\end{definition*}}
\newtheorem{definitions*}[definition*]{Definitions}
\newenvironment{defis*}{\begin{definitions*} \rm}{\end{definitions*}}
\newtheorem{example*}[definition*]{Example}
\newenvironment{exam*}{\begin{example*} \rm}{\end{example*}}
\newtheorem{examples*}[definition*]{Examples}
\newenvironment{exams*}{\begin{examples*} \rm}{\end{examples*}}
\newtheorem{nothing*}[definition*]{$\!\!$}
\newenvironment{noth*}{\begin{nothing*} \rm}{\end{nothing*}}

\theoremstyle{plain}
\newtheorem{prop*}[definition*]{Proposition}
\newtheorem{lemm*}[definition*]{Lemma}
\newtheorem{coro*}[definition*]{Corollary}
\newtheorem{theo*}[definition*]{Theorem}

\theoremstyle{remark}
\newtheorem{remark*}[definition*]{Remark}
\newenvironment{rema*}{\begin{remark*} \rm}{\end{remark*}}
\newtheorem{remarks*}[definition*]{Remarks}
\newenvironment{remas*}{\begin{remarks*} \rm}{\end{remarks*}}

\begin{document}

\def \n {{\underline{n}}}
\def \chara {{{{\rm Char}}}}
\def \ker {{{\rm Ker}}}
\def \Y {{{\mathfrak{Y}}}}
\def \u {{{\mathfrak{u}}}}
\def \v {{{\mathfrak{v}}}}
\def \w {{{\mathfrak{w}}}}
\def \T {{{\mathfrak{T}}}}
\def \Dt {{{\widetilde{D}}}}
\def \im {{{\rm Im}}}
\def \Hom {{{\rm Hom}}}
\def \Cb {{\overline{C}}}
\def \K {{\mathbb{K}}}
\def \Nt {{\widetilde{\cN}}}
\def \Ct {{\widetilde{C}}}
\def \Ch {{\widehat{C}}}
\def \Xh {{\widehat{X}}}
\def \Yh {{\widehat{Y}}}
\def \pt {{\widetilde{p}}}
\def \qt {{\widetilde{q}}}
\def \gl {{\mathfrak{gl}}}
\def \Ft {\widetilde{F}}
\def \Tt {\widetilde{T}}
\def \ltt {\widetilde{\lt}}
\def \Tb {\overline{T}}
\def \ltb {\overline{\lt}}
\def \codim {{\rm Codim}}
\def \sca #1#2{\langle #1,#2 \rangle}
\def\x {{\underline{x}}}
\def\y {{\underline{y}}}
\def\aut{{\rm Aut}}
\def\ra{\rightarrow}
\def\s{\sigma}\def\OO{\mathbb O}
\def\QQ{\mathbb Q}
 \def\CC{\mathbb C} \def\ZZ{\mathbb Z}\def\JO{{\mathcal J}_3(\OO)}
\newcommand{\G}{\mathbb{G}}
\def\proof{\noindent {\it Proof.}\;}
\def\qed{\hfill $\square$}
\def \uh {{\widehat{u}}}
\def \vh {{\widehat{v}}}
\def \fh {{\widehat{f}}}
\def \wh {{\widehat{w}}}
\def \Wh {{{W_{{\rm aff}}}}}
\def \Wt {{\widetilde{W}_{{\rm aff}}}}
\def \Qt {{\widetilde{Q}}}
\def \Waff {{W_{{\rm aff}}}}
\def \Waffm {{W_{{\rm aff}}^-}}
\def \Wpaff {{{(W^P)}_{{\rm aff}}}}
\def \Wtpaff {{{(\widetilde{W}^P)}_{{\rm aff}}}}
\def \Wtaffm {{\widetilde{W}_{{\rm aff}}^-}}
\def \lh {{\widehat{\lambda}}}
\def \pit {{\widetilde{\pi}}}
\def \lt {{{\lambda}}}
\def \xh {{\widehat{x}}}
\def \yh {{\widehat{y}}}


\newcommand{\expxy}{\exp_{x \to y}}
\newcommand{\drat}{d_{\rm rat}}
\newcommand{\dmax}{d_{\rm max}}
\newcommand{\zl}{Z(x,L_x,y,L_y)}


\newcommand{\N}{\mathbb{N}}
\newcommand{\A}{{\mathbb{A}_{\rm Aff}}}
\newcommand{\Ah}{{\mathbb{A}_{\rm Aff}}}
\newcommand{\At}{{\widetilde{\mathbb{A}}_{\rm Aff}}}
\newcommand{\Ht}{{{H}^T_*(\Omega K^{\ad})}}
\renewcommand{\H}{{{H}^T_*(\Omega K)}}
\newcommand{\Ih}{{I_{\rm Aff}}}
\newcommand{\psit}{{\widetilde{\psi}}}
\newcommand{\xit}{{\widetilde{\xi}}}
\newcommand{\Jt}{{\widetilde{J}}}
\newcommand{\Zt}{{\widetilde{Z}}}
\newcommand{\Xt}{{\widetilde{X}}}
\newcommand{\Yt}{{\widetilde{Y}}}
\newcommand{\at}{{\widetilde{A}}}
\newcommand{\Z}{\mathbb Z}
\newcommand{\R}{\mathbb{R}}
\newcommand{\Q}{\mathbb{Q}}
\newcommand{\C}{\mathbb{C}}
\newcommand{\F}{\mathbb{F}}
\newcommand{\p}{\mathbb{P}}
\newcommand{\co}{{\cal O}}

\renewcommand{\a}{{\alpha}}
\newcommand{\az}{\a_\Z}
\newcommand{\ak}{\a_k}

\newcommand{\rc}{\R_\C}
\newcommand{\cc}{\C_\C}
\newcommand{\hc}{\H_\C}
\newcommand{\oc}{\O_\C}

\newcommand{\rk}{\R_k}
\newcommand{\ck}{\C_k}
\newcommand{\hk}{\H_k}
\newcommand{\ok}{\O_k}

\newcommand{\rz}{\R_Z}
\newcommand{\cz}{\C_Z}
\newcommand{\hz}{\H_Z}
\newcommand{\oz}{\O_Z}

\newcommand{\RR}{\R_R}
\newcommand{\CR}{\C_R}
\newcommand{\HR}{\H_R}
\newcommand{\OR}{\O_R}

\newcommand{\re}{\mathtt{Re}}

\newcommand{\matttr}[9]{
\left (
\begin{array}{ccc}
{} \hspace{-.2cm} #1 & {} \hspace{-.2cm} #2 & {} \hspace{-.2cm} #3 \\
{} \hspace{-.2cm} #4 & {} \hspace{-.2cm} #5 & {} \hspace{-.2cm} #6 \\
{} \hspace{-.2cm} #7 & {} \hspace{-.2cm} #8 & {} \hspace{-.2cm} #9
\end{array}
\hspace{-.15cm}
\right )   }


\newcommand{\dual}{{\bf v}}
\newcommand{\com}{\mathtt{Com}}
\newcommand{\rg}{\mathtt{rg}}
\newcommand{\pu}{{\mathbb{P}^1}}
\newcommand{\scal}[1]{\langle #1 \rangle}
\newcommand{\MK}[2]{{\overline{{\rm M}}_{#1}(#2)}}
\newcommand{\mor}[2]{{{\rm Mor}_{#1}(\pu,#2)}}

\newcommand{\fg}{\mathfrak g}
\newcommand{\fgad}{{\mathfrak g}^{\rm ad}}
\renewcommand{\fh}{\mathfrak h}
\newcommand{\fu}{\mathfrak u}
\newcommand{\fz}{\mathfrak z}
\newcommand{\fn}{\mathfrak n}
\newcommand{\fe}{\mathfrak e}
\newcommand{\fp}{\mathfrak p}
\newcommand{\ft}{\mathfrak t}
\newcommand{\fl}{\rm Fl}
\newcommand{\fq}{\mathfrak q}
\newcommand{\fsl}{\mathfrak {sl}}
\newcommand{\fgl}{\mathfrak {gl}}
\newcommand{\fso}{\mathfrak {so}}
\newcommand{\fsp}{\mathfrak {sp}}
\newcommand{\ff}{\mathfrak {f}}

\newcommand{\ad}{{\rm ad}}
\newcommand{\jad}{{j^\ad}}
\newcommand{\id}{{\rm id}}


\newcommand{\dynkinadeux}[2]
{
$#1$
\setlength{\unitlength}{1.2pt}
\hspace{-3mm}
\begin{picture}(12,3)
\put(0,3){\line(1,0){10}}
\end{picture}
\hspace{-2.4mm}
$#2$
}

\newcommand{\mdynkinadeux}[2]
{
\mbox{\dynkinadeux{#1}{#2}}
}

\newcommand{\dynkingdeux}[2]
{
$#1$
\setlength{\unitlength}{1.2pt}
\hspace{-3mm}
\begin{picture}(12,3)
\put(1,.8){$<$}
\multiput(0,1.5)(0,1.5){3}{\line(1,0){10}}
\end{picture}
\hspace{-2.4mm}
$#2$
}

\newcommand{\poidsesix}[6]
{
\hspace{-.12cm}
\left (
\begin{array}{ccccc}
{} \hspace{-.2cm} #1 & {} \hspace{-.3cm} #2 & {} \hspace{-.3cm} #3 &
{} \hspace{-.3cm} #4 & {} \hspace{-.3cm} #5 \vspace{-.13cm}\\
\hspace{-.2cm} & \hspace{-.3cm} & {} \hspace{-.3cm} #6 &
{} \hspace{-.3cm} & {} \hspace{-.3cm}
\end{array}
\hspace{-.2cm}
\right )      }

\newcommand{\copoidsesix}[6]{
\hspace{-.12cm}
\left |
\begin{array}{ccccc}
{} \hspace{-.2cm} #1 & {} \hspace{-.3cm} #2 & {} \hspace{-.3cm} #3 &
{} \hspace{-.3cm} #4 & {} \hspace{-.3cm} #5 \vspace{-.13cm}\\
\hspace{-.2cm} & \hspace{-.3cm} & {} \hspace{-.3cm} #6 &
{} \hspace{-.3cm} & {} \hspace{-.3cm}
\end{array}
\hspace{-.2cm}
\right |      }

\newcommand{\poidsesept}[7]{
\hspace{-.12cm}
\left (
\begin{array}{cccccc}
{} \hspace{-.2cm} #1 & {} \hspace{-.3cm} #2 & {} \hspace{-.3cm} #3 &
{} \hspace{-.3cm} #4 & {} \hspace{-.3cm} #5 & {} \hspace{-.3cm} #6
\vspace{-.13cm}\\
\hspace{-.2cm} & \hspace{-.3cm} & {} \hspace{-.3cm} #7 &
{} \hspace{-.3cm} & {} \hspace{-.3cm}
\end{array}
\hspace{-.2cm}
\right )      }

\newcommand{\copoidsesept}[7]{
\hspace{-.12cm}
\left |
\begin{array}{cccccc}
{} \hspace{-.2cm} #1 & {} \hspace{-.3cm} #2 & {} \hspace{-.3cm} #3 &
{} \hspace{-.3cm} #4 & {} \hspace{-.3cm} #5 & {} \hspace{-.3cm} #6
\vspace{-.13cm}\\
\hspace{-.2cm} & \hspace{-.3cm} & {} \hspace{-.3cm} #7 &
{} \hspace{-.3cm} & {} \hspace{-.3cm}
\end{array}
\hspace{-.2cm}
\right |      }

\newcommand{\poidsehuit}[8]{
\hspace{-.12cm}
\left (
\begin{array}{cccccc}
{} \hspace{-.2cm} #1 & {} \hspace{-.3cm} #2 & {} \hspace{-.3cm} #3 &
{} \hspace{-.3cm} #4 & {} \hspace{-.3cm} #5 & {} \hspace{-.3cm} #6 &
{} \hspace{-.3cm} #7   \vspace{-.13cm}\\
\hspace{-.2cm} & \hspace{-.3cm} & {} \hspace{-.3cm} #8 &
{} \hspace{-.3cm} & {} \hspace{-.3cm}
\end{array}
\hspace{-.2cm}
\right )      }

\newcommand{\copoidsehuit}[8]{
\hspace{-.12cm}
\left |
\begin{array}{cccccc}
{} \hspace{-.2cm} #1 & {} \hspace{-.3cm} #2 & {} \hspace{-.3cm} #3 &
{} \hspace{-.3cm} #4 & {} \hspace{-.3cm} #5 & {} \hspace{-.3cm} #6 &
{} \hspace{-.3cm} #7  \vspace{-.13cm}\\
\hspace{-.2cm} & \hspace{-.3cm} & {} \hspace{-.3cm} #8 &
{} \hspace{-.3cm} & {} \hspace{-.3cm}
\end{array}
\hspace{-.2cm}
\right |      }


\def\cA{{\cal A}} \def\cC{{\cal C}} \def\cD{{\cal D}} \def\cE{{\cal E}}
\def\cFt{{\widetilde{{\cal F}}}}
 \def\cG{{\cal G}} \def\cH{{\cal H}} \def\cI{{\cal I}}
\def\cK{{\cal K}} \def\cL{{\cal L}} \def\cM{{\cal M}} \def\cN{{\cal N}}
\def\cO{{\cal O}}
\def\cP{{\cal P}} \def\cQ{{\cal Q}} \def\cT{{\cal T}} \def\cU{{\cal U}}
\def\cV{{\cal V}} \def\cX{{\cal X}} \def\cY{{\cal Y}} \def\cZ{{\cal Z}}


\newcommand{\KK}{\mathbb{K}}
\newcommand{\PP}{\mathbb{P}}
\newcommand{\cF}{\mathcal{F}}


\title{Springer fiber components in the two columns case \\ for types $A$ and $D$ are normal}
\author{Nicolas Perrin and Evgeny Smirnov}

\maketitle

\begin{abstract}
We study the singularities of the irreducible components of the
Springer fiber over a nilpotent element $N$ with $N^2=0$ in a Lie algebra
of type $A$ or $D$ (the so-called two columns case).
We use Frobenius splitting techniques to prove that these irreducible
components are normal, Cohen--Macaulay, and have rational singularities.
\end{abstract}

 {\def\thefootnote{\relax}
 \footnote{ \hspace{-6.8mm}
 }
 }

\section{Introduction}

Let $\K$ be an algebraically closed field of arbitrary
characteristic not equal to 2. Let $N$ be a nilpotent element in a
Lie algebra $\fg=\gl_n$ (type $A$) or $\fg=\so_{2n}$ (type $D$). We
consider the Springer fiber $\cF_N$ over $N$. It is the fiber of the
famous Springer resolution of the nilpotent cone $\cN\subset\fg$
over $N$.

This resolution can be constructed as follows. Let $\cF$ be the variety
of complete flags in $\K^n$ (resp. $\Or\cF$ the variety of complete
\emph{isotropic} flags, see Section \ref{section-ortho} for the
description of the Springer fiber $\Or\cF_N$ in this case).
A flag $f=(V_i)_{i\in[0,n]}$,
where $V_i$ is a vector subspace of $\K^n$ of dimension $i$, is stabilized 
by $N\in\cN$ if $N(V_i)\subset V_{i-1}$ for all $i>0$. We shall denote this
by $N(f)\subset f$. Define the variety
$$\Nt=\{(f,N)\in\cF\times\cN\ \mid N(f)\subset f\}.$$
The projection $\Nt\to\cF$ is a smooth morphism
thus $\Nt$ is smooth. The natural projection $\Nt\to\cN$ is birational and proper. It is a resolution of singularities for $\cN$ called the
\emph{Springer resolution}.

The \emph{Springer fibers}, i.e., the fibers of the Springer resolution,
are of great interest. They are connected (this can be seen directly
or follows from the normality of the nilpotent cone $\cN$), equidimensional,
but not irreducible. There is a natural combinatorial framework to
describe them: Young diagrams and standard tableaux.

The irreducible components of the Springer fibers are not well
understood. For example, it is known that in general the components
are singular but there is no general description of the singular
components. There are only partial answers in type $A$. First, it is
known in the so-called  \emph{hook} and \emph{two lines} cases that all
the components are smooth (see \cite{fung}). The first case where
singular components appear is the \emph{two columns} case. A
description of the  singular components in the two columns case has
been given by L. Fresse in \cite{fresse} and \cite{fresse_CR}. In
their recent work \cite{fresse2} L. Fresse and A. Melnikov describe
the Young diagrams for which all irreducible components are smooth.

In this paper, we focus on the the two columns case, that is to say, the case
of nilpotent elements $N$ of order 2. The corresponding
Young diagram $\lambda=\lambda(N)$ has two columns.
We want to
understand the type of singularities appearing in a component of the
Springer fiber.

Let $X$ be an irreducible component of the Springer fiber $\cF_N$,
resp. $\Or\cF_N$, in type $A$, resp. $D$, with $N$ nilpotent such that $N^2=0$.
In the two columns case, we describe a resolution $\pi\colon \Xt\to X$ of the
irreducible component $X$. We use this resolution  to prove, for
$\chara(\K)>0$,  that $X$ is Frobenius split, and deduce the following
result for arbitrary characteristic:

\begin{theo}
\label{main1}
The irreducible component $X$ is normal.
\end{theo}

We are able to prove more on the resolution $\pi$. Recall that a proper
birational morphism $f:X\to Y$ is called a rational
  resolution if $X$ is smooth and if the equalities $f_*\co_X=\co_Y$
  and $R^if_*\co_X= R^if_*\omega_X=0$ for $i>0$ are satisfied.
We prove the following

\begin{theo}
\label{main2}
The morphism $\pi$ is a rational resolution.
\end{theo}

\begin{coro}
\label{main3}
The irreducible component $X$ is Cohen--Macaulay with dualizing sheaf
$\pi_*\omega_\Xt$.
\end{coro}

Rational singularities are well defined in characteristic zero. In
this case we obtain the following

\begin{coro}
\label{main4}
  If $\chara(\K)=0$, then $X$ has rational singularities.
\end{coro}

\paragraph*{Acknowledgements.} We express our gratitude to Michel
Brion for useful discussions, in particular, concerning the existence of
Frobenius splittings using the pair $({\rm SL}_{2n},{\rm
  Sp}_{2n})$. We thank X. He and J.F. Thomsen for giving us a preliminary version of their work \cite{he-thomsen} which simplifies the proof of Theorem \ref{theo-split}. We also thank Catharina Stroppel for discussions on her
paper \cite{str_web} which were the starting point of this
project. Finally, we are grateful to the referee for his careful reading
of our paper and for correcting some mistakes in the initial version
of this paper, especially Proposition \ref{prop-descrip}. E.S. was
partially supported by the RFBR grant 10-01-00540 and by the
RFBR--CNRS grant 10-01-93111. 

\tableofcontents

\section{Irreducible components of Springer fibers in type $A$}

\subsection{General case}

Let $N\in\gl_n$ be a nilpotent element, and let $(m_i)_{i\in[1,r]}$ be the
sizes of Jordan blocks of $N$. To $N$ we assign a Young diagram
$\lambda=\lambda(N)$ of size $n$ with rows of lengths $(m_i)_{i\in[1,r]}$.
We refer to \cite{fulton} for more details on Young diagrams.

\begin{defi}
A \emph{standard Young tableau} of shape $\lambda$ is a bijection
$\tau\colon\lambda\to [1,n]$ such that the numbers assigned to the boxes in
each row are decreasing from left to right, and the numbers in each
column are decreasing from top to bottom.
\end{defi}

\begin{rema}
Usually, one requires that the integers in the boxes of a standard
tableau \emph{increase}, not decrease from left to right and from top
to bottom. However, using \emph{decreasing} tableaux in our case
simplifies the notation, so we decided to follow this (rather unusual)
definition.
\end{rema}

Remark that the datum of a standard tableau $\tau$ is equivalent to the
datum of a chain of decreasing Young diagrams
$\lambda=\lambda^{(0)}\supset\lambda^{(1)}\supset\lambda^{(2)}\supset\dots\supset\lambda^{(n)}=\emptyset$, where $\lambda^{(i)}$ is the set of
the $n-i$ boxes with the largest numbers, that is,
$\tau^{-1}(\{i+1,\dots,n\})$.

Let $f=(V_i)\in\cF$ be an $N$-stable flag. We assign to it a
standard tableau of shape $\lambda=\lambda(N)$ in the following way.
Consider the quotient spaces $V^{(i)}=V/V_i$. The endomorphism $N$
induces an endomorphism of each of these quotients $N^{(i)}\colon
V^{(i)}\to V^{(i)}$. Take the Young diagram $\lambda^{(i)}$
corresponding to $N^{(i)}$; it consists of $n-i$ boxes. Clearly,
$\lambda^{(i)}$ differs from $\lambda^{(i-1)}$ by one corner box. So
we obtain a chain of decreasing Young diagrams, which is equivalent
to a standard Young tableau $\tau(f)$.

Let $\tau$ be a standard Young tableau of shape $\lambda(N)$. Define
$$
X_\tau^0=\{f\in\cF_N\mid \tau(f)=\tau\}.
$$

The following theorem is due to Spaltenstein \cite{spal}.

\begin{theo} For each standard tableau $\tau$,  the subset $X^0_\tau$
 is a smooth irreducible subvariety of $\cF_N$.
 Moreover, $\dim X^0_\tau=\dim\cF_N$, so $X_\tau=\overline{X^0_\tau}$
 is an irreducible component of $\cF_N$. Any irreducible component of $\cF_N$
is obtained in this way.
\end{theo}

\subsection{Two columns case}

In this paper, we focus on the case of nilpotent elements $N$
such that $N^2=0$. This is equivalent to saying that the Young
diagram $\lambda(N)$ consists of (at most) two columns. Denote by $r$ the rank
of $N$ or, equivalently, the number of boxes in the second column.
Let $X=X_\tau$ be the irreducible component of the Springer fiber
over $N$ corresponding to a standard tableau $\tau$. Denote the
increasing sequence of labels in the second column of the standard
tableau $\tau$ by $(p_i)_{i\in[1,r]}$. Set $p_0=0$ and
$p_{r+1}=n+1$.

According to F.Y.C. Fung \cite{fung}, the previous Theorem can be
reformulated as follows.

\begin{prop}
The irreducible component $X$ is the closure of the variety
$$X^0=\left\{(V_i)_{i\in[0,n]}\in\cF_N\ \left|
\begin{array}{ll}
V_{i}\subset V_{i-1}+\im N&\textrm{for $i\in\{p_1,\cdots,p_{r}\}$}
\\
V_i\not\subset V_{i-1}+\im N&\textrm{otherwise}\\
\end{array} \right.\right\}.$$
\end{prop}

An easy interpretation of this result is the following

\begin{coro}
\label{coro-comp}
The irreducible component $X$ is the closure of the variety
$$X^0=\left\{(V_i)_{i\in[0,n]}\in\cF_N\ \mid
\dim(\im N\cap V_{i})=k \textrm{ for all $k\in[0,r]$ and all
$i\in[p_k,p_{k+1})$}\right\}.$$
\end{coro}

\begin{proo}
  We prove this by induction on $i$. We have $\dim(\im N\cap
  V_0)=0$. The result is implied by the following equivalence:
  $(V_{i+1}\subset V_i+\im N) \Leftrightarrow (\dim(\im N\cap V_{i+1}) =
\dim(\im N\cap V_i)+1)$.
\end{proo}

\subsection{A birational transformation of the Springer fiber}
\label{transf-1}

The above description gives a natural way to construct a resolution of
singularities for $X$. We start with the following simple birational
transformation of $X$. Define the variety $\Xh$ as follows:
$$
\Xh=\{((F_k)_{k\in[0,r]},(V_i)_{i\in[0,n]})\in\cF(\im N)\times\cF\
\mid F_k\subset V_{p_k}\subset N^{-1}(F_{k-1}),\ \forall
k\in[1,r]\},
$$
where $\cF(\im N)$ denotes the variety of complete flags in $\im N$.
The natural projections of the product $\cF(\im N)\times\cF$ on its
two factors induce two maps $p_X:\Xh\to\cF$ and $q_X:\Xh\to\cF(\im N)$.

One of the
main features of the two columns case that we will use is the following easy
observation: $\im N\subset\ker N$. In particular, for any flag
$(F_k)_{k\in[0,r]}\in\cF(\im N)$, the equalities $F_r=\im N$ and
$N^{-1}(F_0)=\ker N$ imply the following inclusions:
$$
F_0\subset\cdots\subset F_r\subset N^{-1}(F_0)\subset\cdots\subset
N^{-1}(F_r).
$$
Fixing subspaces $(F_i)_{i\in[r,n-r]}$ with $\dim(F_i)=i$ such that
$$\im N\subset F_r\subset\cdots\subset F_{n-r}\subset\ker N$$
gives for any choice of $(F_k)_{k\in[0,r]}\in\cF(\im N)$ a complete flag
$$F_0\subset\cdots\subset F_r\subset F_{r+1}\subset\cdots\subset
F_{n-r-1}\subset N^{-1}(F_0)\subset\cdots\subset N^{-1}(F_r)$$
in $N^{-1}(F_r)=\K^n$. We denote this complete flag by $F_\bullet$.

\begin{prop}
\label{prop-p-birat}
(\i) The map $q_X$ is dominant and is a locally trivial fibration over $\cF(\im N)$. Its
fiber over $(F_k)_{k\in[0,r]}$ is isomorphic to the following Schubert
variety associated to $F_\bullet$:
$$
\cF_w=\{(V_i)_{i\in[0,n]}\in\cF\ \mid F_k\subset V_{p_k}\subset
N^{-1}(F_{k-1}),\ \forall k\in[1,r]\}.$$

(\i\i) The map $p_X$ is birational onto $X$.
\end{prop}

\begin{proo}
\textit{(\i)} The first part is clear from the definition of $\Xh$.

\textit{(\i\i)} Let $(V_i)_{i\in[0,n]}$ be in $X^0$. We may define $F_k=\im
N\cap V_{p_k}$ for $k\in[0,r]$. We have $\dim F_k=k$. Since
$(V_i)_{i\in[0.n]}$ is
in the Springer fiber, we also have the inclusion $N(V_{p_k})\subset
V_{p_k-1}$. But $N(V_{p_k})\subset\im N$, thus
$$N(V_{p_k})\subset\im N\cap V_{p_k-1}.$$
Since $(V_i)_{i\in[0,n]}$ is in $X^0$, we have $\im N\cap V_{p_k-1}=\im
N\cap V_{p_{k-1}}$. Therefore we have the inclusion:
$$V_{p_k}\subset N^{-1}(\im N \cap V_{p_{k-1}})=N^{-1}(F_{k-1}).$$
In particular $X^0$ is contained in the image of $p_X$.

Conversely, let $(F_k)_{k\in[0,r]}\in\cF(\im N)$ and
$(V_i)_{i\in[0,n]}$ in the Schubert variety $\cF_w$ associated to
$F_\bullet$. It is easy to check that for $(V_i)_{i\in[0,n]}$
general in the Schubert variety, we have $\im N\cap V_i=F_k$ for
$i\in[p_k,p_{k+1})$. Furthermore, for $i\in[p_k,p_{k+1})$ we have
the inclusions
$$N(V_{i+1})\subset N(V_{p_{k+1}})\subset \im N\cap V_{p_k}= F_{k}
\subset V_{p_k}\subset V_i,$$
therefore $(V_i)_{i\in[0,n]}$ is in $X^0$.
\end{proo}

\subsection{A Schubert variety containing $X$}
\label{schub}

Let us consider the following subvariety of $\cF(\im N)\times\cF$
containing $\Xh$:
$$\Yh=\{((F_k)_{k\in[0,r]},(V_i)_{i\in[0,n]})\in\cF(\im N)\times\cF\ \mid
F_k\subset V_{p_k},\ \forall
k\in[1,r]\}.$$
As for $\Xh$, the natural projections of the product $\cF(\im
N)\times\cF$ on its two factors induce two maps $p_Y:\Yh\to\cF$
and $q_Y:\Yh\to\cF(\im N)$.

\begin{prop}
(\i) The map $q_Y$ is dominant and is a locally trivial fibration with
fiber over $(F_k)_{k\in[0,r]}$ isomorphic to the following Schubert
variety associated to $F_\bullet$:
$$\cF_v=\{(V_i)_{i\in[0,n]}\in\cF\ \mid
F_k\subset V_{p_k},\ \forall
k\in[1,r]\}.$$

(\i\i) The map $p_Y$ is birational onto the Schubert variety
$$Y=\{(V_i)_{i\in[0,n]}\in\cF\ \mid
\dim(\im N\cap V_{p_k})\geq k,\ \forall
k\in[1,r]\}.$$
\end{prop}

\begin{proo}
\textit{(\i)} The first part is clear from the definition of $\Yh$.

\textit{(\i\i)} The image of $p_Y$ is contained in $Y$. Conversely, let
$(V_i)_{i\in[0,n]}$ be general in $Y$. We then have $\dim(\im N \cap
V_{p_k})=k$ and we may define $F_k=\im N\cap V_{p_k}$ for
$k\in[0,r]$. We have $\dim F_k=k$ and
$((F_k)_{k\in[0,r]},(V_i)_{i\in[0,n]})$ is in the fiber of $p_Y$ over
$(V_i)_{i\in[0,n]}$.
\end{proo}

\section{Irreducible components of Springer fibers in type $D$}
\label{section-ortho}

\subsection{Preliminaries on orthogonal groups and Springer fibers}

Let $V$ be a $2n$-dimensional vector space. Consider the group
$\SO(V)$ of unimodular linear operators  preserving a symmetric
nondegenerate bilinear form $\omega$. Let $B$ be a Borel subgroup in
$\SO(V)$. The flag variety $\SO(V)/B$ is the variety $\Or\cF$ of
\emph{orthogonal flags} defined by
$$
\Or\cF=
\{V_0\subset V_1\subset\dots\subset V_{n-1}\subset
V_{n+1}\dots\subset V_{2n}\mid
V_{2n-i}=V_i^\perp \textrm{ and } \dim V_i=i\ \textrm{for }i\leq n-1\}.$$
We will consider elements in $\Or\cF$ as $n$-tuples of
nested isotropic vector spaces $((V_i)_{i\in[0,n-1]})$. We
recover the usual notion of orthogonal flags because there are
exactly two maximal isotropic subspaces between $V_{n-1}$ and its
orthogonal $V_{n+1}=(V_{n-1})^\perp$.

Let $N\in\gl(V)$ be a nilpotent
element. $N$ is said to be \emph{orthogonalizable} if there exists a
symmetric nondegenerate bilinear form $\omega$ on $V$ such that $N$
is $\omega$-invariant; that is,
$$
\omega(Nv,w)+\omega(v,Nw)=0.
$$
This means that $N\in\so(V)$, where $\so(V)$ is the set of elements
of $\gl(V)$ leaving $\omega$ invariant.
The following easy consequence of the Jacobson--Morozov Theorem can
be found, for instance, in \cite[Chap.~6, 2.3]{viniti90}.

\begin{prop}\label{nilpotent_ortho} A nilpotent element $N$ is
  orthogonalizable if in
the corresponding partition $Y(N)$ each even term occurs with even
multiplicity (such partitions will be called \emph{admissible}).
\end{prop}

\begin{definition} Given a nilpotent element $N\in\so(V)$, we define a
  \emph{Springer
fiber of type $D$} in the usual way: namely, as the set of all
orthogonal flags stabilized by $N$:
$$\Or\cF_N=\{((V_i)_{i\in[0,n-1]})\in\Or\cF\mid N(V_i)\subset
V_{i-1} \text{ for } i\in[0,n-1] \textrm{ and } N(V_{n-1}^\perp)\subset
V_{n-1}\}.$$
\end{definition}

Remark that the orthogonalizability condition on $N$ implies that
$N(V_i^\perp)\subset V_{i+1}^\perp$ for $i\leq n-2$. A description of irreducible
components of Springer fibers in types $B$, $C$, and $D$ was given in
M.~van~Leeuwen's Ph.D.~thesis \cite{leeuwen}. We briefly recall this
description here for the type $D$.

\begin{definition}\label{SDT} Let $\lambda$ be a Young diagram with $2n$ boxes. A map $\tau$ from the boxes of $\lambda$
to $[1,n]$ is called a \emph{standard domino tableau}, if the
following conditions hold:

(\i) For each $i$, the pre-image $\tau^{-1}(i)$ consists exactly of
two adjacent boxes (adjacent either by horizontal or by vertical);

(\i\i) For each $i$, the set of boxes
$\lambda^{(i)}(N):=\tau^{-1}([i+1,n])$ corresponding to the numbers
greater than $i$ forms a Young diagram.

Moreover, a standard domino tableau is said to be \emph{admissible},
if all the diagrams $\lambda^{(i)}(N)$ are admissible (in the sense of
Prop.~\ref{nilpotent_ortho}).
\end{definition}

We will think of the pair of boxes $\tau^{-1}(i)$ as of a domino
tile indexed by the number $i$. Each of these tiles can be either horizontal or vertical.

\begin{example}
Let $\lambda=(3,3)$. Then there are three standard domino tableaux of
shape $\lambda$ (see below), but only the first two of them are
admissible. Indeed, for the third diagram $\tau^{-1}(3)$ corresponds
to the Young diagram with one row of length 2, which is not
admissible.
$$
\young(321,321)  \qquad \young(322,311) \qquad \young(331,221)
$$
\end{example}

\begin{definition}
Let $\tau$ be an admissible standard domino tableau of shape
$\lambda(N)$. We assign to it a subset $X_\tau$ of the Springer
fiber $\Or\cF_N$ obtained as the closure of the set $X_\tau^0$ of
flags $(V_i)_{i\in[0,n-1]}$ in $\Or\cF_N$ such that
$N|{}_{V_i^\perp/V_i}$ corresponds to the partition
$\tau^{-1}([i+1,n])$ for each $i<n$.

By definition of $\Or\cF_N$, $N$ is well defined on
$V_i^\perp/V_i=V_{2n-i}/V_i$, so this makes sense.
\end{definition}

The following theorem is due to M.~van Leeuwen \cite[Sec.~3.2]{leeuwen}.

\begin{theo}\label{domino}
$X_\tau$ is an irreducible component of $\Or\cF_N$; all its irreducible
components are obtained in this way. In particular, there is a
bijection between the admissible standard domino tableaux of shape
$\lambda(N)$ and the irreducible components of $\Or\cF_N$.
\end{theo}

\subsection{Description of components in the two columns case}

Throughout this subsection we fix a nilpotent element $N\in\so(V)$
such that $N^2=0$, and an admissible standard domino tableau $\tau$ of shape
$\lambda(N)$. We begin with the following combinatorial observation.

\begin{prop}\label{vertical_tiles} 
The Young diagram $\lambda(N)$ has at most two columns. Each admissible standard domino tableau of shape $\lambda(N)$ contains only vertical tiles.
\end{prop}

Note that in particular, the rank of $N$ has to be even. Let ${\rm rk}(N)=2r$, and let $(p_i)_{i\in[1,r]}$ be the numbers of domino tiles forming the second column of the diagram $Y(N)$. We formally set $p_0=0$, $p_{r+1}=2n+1$.

Now we endow the subspace $\im N$ with a bilinear
form $\alpha$ as follows. For $u,v\in\im N$,
$$
\alpha(u,v)=\omega(u,v'),\text{ where }v'\in N^{-1}(v).
$$

\begin{prop} 
The form $\alpha$ is a skew-symmetric nondegenerate form on $\im N$.
\end{prop}

\begin{proo}
  We readily see that $\alpha$ is well-defined. To show that it
is skew-symmetric, take two vectors $u,v\in\im N$ along with their
preimages $u'\in N^{-1}(u)$ ,$v'\in N^{-1}(v)$. Then
$$
\alpha(u,v)=\omega(N(u'),v')=-\omega(u',N(v'))=-\omega(N(v'),u')=-\alpha(v,u).
$$
The non-degeneracy of $\alpha$ is also obvious.
\end{proo}

\begin{rema} This is a particular case of the construction of
a family of nondegenerate bilinear forms on $(\Ker N\cap \im
N^i)/(\Ker N\cap \im N^{i+1})$, which works for arbitrary nilpotent
$N\in\so(V)$. See \cite[Section~2.3]{leeuwen} for details.
\end{rema}

We shall denote by $\angle$ the orthogonality relation for the form
$\alpha$. We consider the \emph{symplectic flag variety}  $\Sp\cF(\im
N)$, defined as follows:
$$
\Sp\cF(\im N)=\{(0=F_0\subset F_1\subset F_2\subset \dots\subset
F_{2r}=\im N)\mid F_{2r-k}=F_k^\angle\}.
$$
Similarly as for $\Or\cF$, an element of $\Sp\cF$ can also be seen as a sequence 
$(F_k)_{k\in[0,r]}$ of nested $\alpha$-isotropic subspaces.

As in the type $A$ case (see Corollary \ref{coro-comp}), the description of
irreducible components can be reformulated as follows. 

\begin{prop}
\label{prop-descrip}
(\i) Let $f=(V_i)_{i\in[0,n-1]}\in\Or\cF_N$ be an orthogonal $N$-stable flag. There is a unique partial flag $(U_i)_{i\in[0,n-1]}$ of $\im N$ such that for all $i$, we have $U_i\subset V_i\cap N(U_i^\perp)$, with $U_i$ maximal with this property. Moreover, $U_i$ is given by $U_i=V_i\cap N(U_{i-1}^\perp)$
for all $i$.

(\i\i) The subspaces $U_i$ are $\alpha$-isotropic. 

(\i\i\i) $f\in X_\tau^0$ if and only if $\dim U_i=\#(\{p_1,\dots,p_r\}\cap[1,i])$. In that case we have the equality $U_i=\sum_{j=0}^{i}V_j\cap N(V_j^\perp)$.
\end{prop}

Note that the indices $i$ do not correspond to the dimension of $U_i$ and that some of these subspaces may coincide.

\begin{proo}
 \textit{(\i)} Let us begin with the following lemma.

\begin{lemm}
\label{fact-ref}
(\i) A subspace $W\subset V_i$ is maximal for the property $W\subset N(W^\perp)$ if and only if it satisfies $W_i=V_i\cap N(W^\perp)$. 

(\i\i) If $W'\subset W$ are maximal such that $W'\subset V_{i-1}\cap N(W'^\perp)$ and $W\subset V_{i}\cap N(W^\perp)$, then $\dim W\in\{\dim W',\dim W'+1\}$.
\end{lemm}

\begin{proo}
\textit{(\i)} If $W$ is maximal, let $U=V_i\cap N(W^\perp)$, we have $W\subset U$ thus $U\subset V_i\cap N(U^\perp)$ and $W=U$ by maximality. Conversely, if $W=V_i\cap N(W^\perp)$ and if we have $W\subset U$ with $U\subset V_i\cap N(U^\perp)$, then we have the inclusions $W\subset U\subset V_i\cap N(U^\perp)\subset V_i\cap N(W^\perp)$. We thus have equalities for all these inclusions.

\textit{(\i\i)} By \textit{(\i)}, we have the equalities $W'=V_{i-1}\cap N({W'}^\perp)$ and $W=V_i\cap N(W^\perp)$. We thus have the inclusion $W\subset V_i\cap N({W'}^\perp)$ and $W'=V_{i-1}\cap N({W'}^\perp)$ is of codimension at most one in $V_i\cap N({W'}^\perp)$. Thus its codimension in $W$ is also at most one.
\end{proo}

We define the subspaces $U_i$ inductively by the rule $U_0=0$ and $U_i=V_i\cap N(U_{i-1}^\perp)$. Let us check by induction that they are maximal so suppose that $U_{i-1}=V_{i-1}\cap N(U_{i-1}^\perp)$. Either $U_i=U_{i-1}$, and the maximality of $U_i=V_i\cap N(U_i^\perp)$ follows from Lemma~\ref{fact-ref} \textit{(\i)}. Or $U_i= U_{i-1}\oplus\langle x\rangle$. Let $y\in U_{i-1}^\perp$ with $N(y)=x$. We have the equalities $\omega(x,y)=\omega(N(y),y)=-\omega(y,N(y))=-\omega(N(y),y)=-\omega(x,y)$ thus $\omega(x,y)=0$. So $y\in U_i^\perp$ and we have the inclusion $U_i\subset N(U_i^\perp)$, and $U_i$ is maximal by Lemma~\ref{fact-ref} \textit{(\i\i)}. The existence of such a flag is shown.

Let $(U'_i)_{i\in[0,n-1]}$ be another partial flag satisfying the same property. We prove by induction that they coincide. Assume that $U'_{i-1}=U_{i-1}$, we see that $U_i'\subset V_i\cap N(U_i'^\perp)\subset V_i\cap N(U_{i-1}'^\perp)=V_i\cap N(U_{i-1}^\perp)=U_i$. Whence the uniqueness by maximality.

\textit{(\i\i)} We have the inclusion $U_i\subset V_i\cap N(U_i^\perp)$, so in particular $U_i\subset N(U_i^\perp)=U_i^\angle$, and all $U_i$ are $\alpha$-isotropic. Note in particular that $\dim U_i\leq r$ for all $i$.

\textit{(\i\i\i)} Let us first prove one more lemma.

\begin{lemm}
Let $(V_i)_{i\in[0,n-1]}\in X_\tau^0$, we have the equalities: 
$$\dim \sum_{j=0}^{p_k} V_j\cap N(V_j^\perp)= \dim \sum_{j=0}^{p_k-1}
(V_j\cap N(V_j^\perp))+1 \textrm{ for all $k\in [1,r]$ and }$$
$$\sum_{j=0}^{k} V_j\cap N(V_j^\perp)= \sum_{j=0}^{k-1}
(V_j\cap N(V_j^\perp)) \textrm{ for all $k\not\in\{p_1,\cdots,p_r\}$.}$$
\end{lemm}

\begin{proo}
  Let us consider the following inclusions:
$$\xymatrix{N(V_j^\perp)\ \ar@{^(->}[rr]&&N(V_{j-1}^\perp)\\
V_j\cap N(V_j^\perp)\ar@{^(->}[u]&&V_{j-1}\cap
N(V_{j-1}^\perp)\ar@{^(->}[u]\\
&V_{j-1}\cap N(V_j^\perp).\ar@{_(->}[ul]\ar@{^(->}[ur]&\\}
$$
If $j$ is in the interval $[p_k,p_{k+1})$, then the left vertical
inclusion is a codimension $2r-2k$ inclusion. The right vertical
inclusion is of the same codimension except for $j=p_k$ in which case
it is of codimension $2r-2k+2$. All the other inclusions are of
codimension at most one.

Assume first that $j=p_k$, then $V_{p_k}\cap N(V_{p_k}^\perp)$ has to
be of dimension one more than $V_{p_k-1}\cap N(V_{p_k-1}^\perp)$ thus
we have $V_{p_k-1}\cap N(V_{p_k}^\perp)=V_{p_k-1}\cap
N(V_{p_k-1}^\perp)$ and $V_{p_k}\cap N(V_{p_k}^\perp)=V_{p_k-1}\cap
N(V_{p_k-1}^\perp)\oplus\langle x\rangle$ with $x$ in $V_{p_k}$ but not in
$V_{p_k-1}$. Therefore we have
$$\sum_{j=0}^{p_k} V_j\cap N(V_j^\perp)= \sum_{j=0}^{p_k-1}
(V_j\cap N(V_j^\perp))\oplus\langle x\rangle$$
and the result follows in this case.
If $j=k\not\in\{p_1,\cdots,p_r\}$, we distinguish between the two cases. If the top
horizontal inclusion is an equality, then we also have the equality
$V_{k-1}\cap N(V_{k}^\perp)= V_{k-1}\cap N(V_{k-1}^\perp)$ and by
dimension count $V_{k}\cap N(V_{k}^\perp)=V_{k-1}\cap
N(V_{k-1}^\perp)$. If not, then $V_k\cap N(V_k^\perp)$ has to be of
dimension one less than $V_{k-1}\cap N(V_{k-1}^\perp)$ thus
we have $V_{k}\cap N(V_{k}^\perp)=V_{k-1}\cap N(V_{k}^\perp)\subset
V_{k-1}\cap N(V_{k-1}^\perp)$. In any case the result follows.
\end{proo}

Let us finish the proof. Since $V_j\cap N(V_j^\perp)\subset V_j\cap
N(U_j^\perp)=U_j$, we get the inclusion $\sum_{j=0}^i (V_j\cap
N(V_j^\perp))\subset U_i$ for all $i$. Note also that for $j=p_k$, we have
$\sum_{j=0}^{p_k} V_j\cap N(V_j^\perp)= \sum_{j=0}^{p_k-1}
(V_j\cap N(V_j^\perp))\oplus\langle x\rangle$ with $x\not\in V_{p_k-1}$ therefore $x\not\in U_{p_k-1}$ and $\dim U_{p_k}>\dim U_{p_k-1}$. We thus have by the previous fact and this observation the following inequalities and equalities:
$$\dim U_i\geq k \textrm{ and } \dim \sum_{j=0}^{i} V_j\cap N(V_j^\perp)=k \textrm{ for $i\in[p_k,p_{k+1})$.}$$
But $U_{n-1}$ being $\a$-isotropic, we have $\dim U_{n-1}\leq r$ therefore we necessarily have the equality $\sum_{j=0}^i V_j\cap N(V_j^\perp)=U_i$, and \textit{(\i\i\i)} also follows.  
\end{proo}

\begin{rema}
\label{rema-L}
Any complete flag $(V_i)_{i\in[0,n-1]}$ in the Springer fiber belongs to some $X^0_\tau$ for some admissible standard tableau $\tau$. In particular, for any such complete flag, we may assign a \emph{Lagrangian} subspace in $\im N$ for the skew form $\alpha$, namely, $U_{n-1}=\sum_{j=0}^{n-1} V_j\cap N(V^\perp_j)$.
\end{rema}

\begin{coro} 
For $f=(V_i)_{i\in[0,n-1]}\in X_\tau^0\subset\Or\cF_N$, the flag $0\subset U_{p_1}\subset\dots\subset U_{p_r}\subset\im N$ belongs to $\Sp\cF(\im N)$.
\end{coro}

\subsection{Birational transformation of the Springer fiber}

In this subsection we construct a birational transformation of a given
irreducible component $X=X_\tau$, analogous to the one described in
Section \ref{transf-1}. We define it as follows:

\begin{equation}\label{xhat_ortho}
\widehat{X}=\{((F_k)_{k\in[0,r]},(V_i)_{i\in[0,n-1]})\subset\Sp\cF(\im
N)\times\Or\cF\mid F_k\subset V_{p_k}\subset
N^{-1}(F_{k-1})\quad\forall k\in[1,r]\}.
\end{equation}
Let us remark that for a flag $(F_k)_{k\in[0,r]}\in\Sp\cF(\im N)$, we
may consider the partial flag
$$
F_0\subset\cdots\subset F_{2r} =\im
N\subset\ker N=N^{-1}(F_0)\subset N^{-1}(F_1)\subset\cdots \subset
N^{-1}(F_{2r}),
$$
where $F_j:=F_{2r-j}^\angle$ for $j>r$, and since $(F_k)_{k\in[0,r]}$ is isotropic for the
form $\alpha$, we have
$F_k^\perp=N^{-1}(F_k^\angle)=N^{-1}(F_{2r-k})$. Therefore the above
partial flag is isotropic for the quadratic form $\omega$ and we may
complete it to an isotropic complete flag. As for the type $A$,
denote the two natural projections by $p_X\colon \widehat
X\to\Sp\cF(\im N)$ and $q_X\colon \widehat X\to\Or\cF$ respectively.

\begin{prop}
\label{prop-birat-D}
(\i) The map $q_X$ is dominant and a locally trivial
  fibration with fiber isomorphic to the following Schubert variety:
$$
\Or\cF_w=\{(V_i)_{i\in[0,n-1]}\in\Or\cF\mid F_k\subset V_{p_k}\subset
N^{-1}(F_{k-1})\quad\forall k\in [1,r]\};
$$
(\i\i) The map $p_X$ is birational onto $X$.
\end{prop}

\begin{proo}
  \textit{(\i)} This is clear from the definition of $\Xh$.

\textit{(\i\i)} Let $f=(V_i)_{i\in[0,n-1]}\in X_\tau^0$. With notation as in Prop.~\ref{prop-descrip}, we may define the subspaces $F_i$ by the condition $F_i=U_{p_i}$. We have $U_{p_k}=V_{p_k}\cap N(U_{p_k}^\perp)$ thus $U_{p_k}\subset V_{p_k}$. Furthermore, we have $N(V_{p_k})\subset V_{p_k-1}$ and because of the inclusions $V_{p_k}\subset V_{p_k}^\perp\subset V_{p_k-1}^\perp\subset U_{p_k-1}^\perp$, we get $N(V_{p_k})\subset N(U_{p_k-1}^\perp)$. Therefore we have the inclusions $N(V_{p_k})\subset V_{p_k-1}\cap N(U_{p_k-1}^\perp)=U_{p_k-1}$ thus $V_{p_k}\subset N^{-1}(U_{p_k-1})=N^{-1}(U_{p_{k-1}})$. The flag $(F_i)_{i\in[0,r]}$ therefore satisfies the conditions~(\ref{xhat_ortho}). So the map $p_X$ is bijective over the dense subset $X^0_\tau\subset X$, and hence it is birational.
\end{proo}

\subsection{A Schubert-like variety containing $X$}
\label{schub-D}

As in type $A$, let us consider the following subvariety of
$\Sp\cF(\im N)\times\Or\cF$ containing $\Xh$:
$$
\Yh=\{((F_k)_{k\in[0,r]},(V_i)_{i\in[0,n-1]})\in\Sp\cF(\im
N)\times\Or\cF\ \mid F_k\subset V_{p_k},\ \forall k\in[1,r]\}.
$$
We shall also consider projection but with a slight modification. Let ${\cal L}_\a(\im N)$ be the variety of all Lagrangian subspaces of $\im N$ for the form $\a$ (\emph{i.e.} isotropic subspaces for the form $\a$ of maximal dimension $r$). We define the projections $q_Y:\Yh\to\Sp\cF(\im N)$ and $p_Y:\Yh\to{\cal L}_\a(\im N)\times\Or\cF$

\begin{prop}
(\i) The map $q_Y$ is dominant and is a locally trivial fibration with
fiber over $(F_k)_{k\in[0,r]}$ isomorphic to the following Schubert
variety
$$\Or\cF_v=\{(V_i)_{i\in[0,n-1]}\in\Or\cF\ \mid
F_k\subset V_{p_k},\ \forall
k\in[1,r]\}.$$

(\i\i) The map $p_Y$ is birational onto the variety
$$Y=\{(L,(V_i)_{i\in[0,n-1]})\in{\cal L}_\a(\im N)\times\Or\cF\ \mid
\dim(L\cap V_{p_k})\geq k,\ \forall
k\in[1,r]\}$$
which is a locally trivial fibration over ${\cal L}_\a(\im N)$ with fiber over $L\in {\cal L}_\a(\im N)$ being the Schubert variety
$$Y_L=\{(V_i)_{i\in[0,n-1]}\in\Or\cF\ \mid
\dim(L\cap V_{p_k})\geq k,\ \forall
k\in[1,r]\}.$$
\end{prop}

\begin{proo}
\textit{(\i)} The first part is clear from the definition of $\Yh$.

\textit{(\i\i)} The last assertion is also clear from the description of $Y$ so we only need to prove that $Y$ is the image of $p_Y$. The image of $p_Y$ is contained in $Y$. Conversely, let
$(L,(V_i)_{i\in[0,n]})$ be general in $Y$. Then we have the equalities $\dim(L\cap V_{p_k})=k$ for all $k\in[0,r]$. We may therefore 
set $F_k=L\cap V_{p_k}$ for $k\in[0,r]$. 
The point 
$((F_k)_{k\in[0,r]},(V_i)_{i\in[0,n-1]})$ is in the fiber of $p_Y$ over
$(L,(V_i)_{i\in[0,n-1]})$.
\end{proo}

\begin{rema}
\label{rema-norm}
The variety $Y$ is normal with rational singularities as a locally trivial fibration with smooth base and fibers which are normal with rational singularities (Schubert varieties).
\end{rema}

\begin{prop}
There is a closed immersion of $X$ in $Y$.
\end{prop}

\begin{proo}
Recall from Remark \ref{rema-L} that for any element $(V_i)_{i\in[0,n-1]}$ in $X_\tau$, the vector space $U_{n-1}=\sum_{j=0}^{n-1}V_j\cap N(V_j^\perp)$ lies in ${\cal L}_\a(\im N)$. Furthermore, for $(V_i)_{i\in[0,n-1]}$ in $X^0_\tau$, we have $\dim(V_{p_k}\cap U_{n-1})=\dim U_{p_k}=k$ for all $k\in[0,r]$. Therefore the map $X\to Y$ defined by $(V_i)_{i\in[0,n-1]}\mapsto((V_i)_{i\in[0,n-1]},L)$ with $L=U_{n-1}$ is a closed immersion.
\end{proo}

\section{Frobenius splitting}

In this section, we assume $\chara(\K)=p>0$, we shall intensively use the
results from the
book \cite{BK}. We refer to this book for the notions of Frobenius
splitting of a scheme $X$ and $B'$-canonical splitting of a scheme $X$
with an action of a Borel subgroup $B'$ of a reductive group $G'$.

We will now deal with type $A$ and type $D$ simultaneously and keep
the notation of the previous two sections. In particular $X$ will
denote a fixed irreducible component of the Springer fiber with two
columns in type $A$ or $D$. To simplify the notation we will denote by
$\cF_w$ both the Schubert variety $\cF_w$ in type $A$ and the Schubert
variety $\Or\cF_w$ in type $D$. We denote by $\n$ the number $n$ in
type $A$ and $n-1$ in type $D$.

\subsection{Bott--Samelson resolutions}
\label{section-BS}

We will use Bott--Samelson resolutions of the Schubert varieties
$\cF_w$ and $\cF_v$ to construct resolutions of $\Xh$ and $\Yh$ and thus
of $X$ and $Y$. Let us fix some notation and recall some basic facts
on Bott-Samelson resolutions (for details we refer to \cite{demazure}
or \cite{BK}).

Recall that the Schubert varieties in $\cF$ are indexed by the elements $u$ of
the Weyl group $W$.
The inclusion of Schubert varieties induces an order on $W$ called the
\emph{Bruhat order}. Any element $u\in W$ can be written as a product
$s_{i_1}\cdots s_{i_k}$ where $s_{i_k}$ is the simple reflection with
respect to the simple root $\a_{i_k}$ (we shall use the notation of
N. Bourbaki here \cite{bou}).
An expression of minimal length is called \emph{reduced},
and its length $k$ is called the \emph{length} of $u$. Let
us also denote, for $\a$ a simple root, by $\G(\a)$, respectively
$\cF(\a)$, the Grassmannian (classical in type $A$, orthogonal in type
$D$) associated to $\a$,
respectively the partial flag variety of all subspaces except those in
$\G(\a)$. Denote by $p_\a$ the projection $\cF\to\cF(\a)$ and by $q_\a$ the projection $\cF\to\G(\a)$. The fiber of $p_\a$ is isomorphic to $\pu$.

Let $F_\bullet$ be a fixed complete flag (classical in type $A$,
orthogonal in type $D$), and let $\u=(\a_{i_1},\cdots, \a_{i_k})$ be a
sequence of simple roots.
We
construct a variety $\cFt_\u$  from these data. For this we consider
the following elementary construction.

\begin{cons}
\label{cons-elementaire}
Having a simple root $\a$, we first
define a variety
$$\cF_{\a}=\{((V_i)_{i\in[0,\n]},W)
\in\cF\times \G(\a)\ \mid W\in
q_\a(p_\a^{-1}(p_\a((V_{i})_{i\in[0,\n]})))\}.$$ There are two natural maps
$\varphi_\a$ and $\psi_\a$ from $\cF_\a$ to $\cF$ defined by
$\varphi_\a((V_i)_{i\in[0,\n]},W)=(V_i)_{i\in[0,\n]}$ and
$\psi_\a((V_i)_{i\in[0,\n]},W)=(p_\a((V_i)_{i\in[0,\n]}),W).$ Remark
that there is a natural section $\s_\a:\cF\to\cF_\a$ of $\varphi_\a$
given by
$(V_i)_{i\in[0,\n]}\mapsto((V_i)_{i\in[0,\n]},q_\a((V_i)_{i\in[0,\n]}))$.

Let $p_Z\colon Z\to\cF$ be a morphism. We define the variety
$Z_\a$ as the fiber product
$$\xymatrix{Z_\a=Z\times_\cF\cF_\a\ar[r]
\ar[d]&Z\ar[d]^{p_Z}\\
\cF_\a\ar[r]^{\varphi_\a}&\cF.}$$
We denote the projection $Z_\a\to Z$ by $f_{Z_\a}$. The section $\s_\a$
induces a section $\s_{Z_\a}$ of $f_{Z_\a}$. We define the map
$p_{Z_\a}:Z_\a\to\cF$ as the composition of the projection
$Z_\a\to \cF_\a$ with $\psi_\a$.
\end{cons}

The Bott--Samelson variety $\cFt_\u$ is constructed from the sequence of
roots $\u$ and the point $F_\bullet$ in $\cF$. Indeed, we set
$Z_0=\{F_\bullet\}$ with the map $p_{Z_0}\colon Z_0\to \cF$ given by the
inclusion and we define $Z_1=(Z_0)_{\a_{i_1}}$ obtained by the elementary
construction from $p_{Z_0}$ and $\a_{i_1}$. We define by induction
$Z_{j+1}=(Z_j)_{\a_{i_{j+1}}}$ obtained by the elementary construction from
$p_{Z_j}$ and $\a_{i_{j+1}}$.
By definition, the Bott--Samelson variety $\cFt_\u$ is $Z_k$. The map
$f_{Z_j}\colon Z_j\to Z_{j-1}$ is a $\pu$-bundle for all $j$ and therefore
$\cFt_\u$ is smooth. The sections $\s_{Z_j}$ define divisors
$D_j=f_{Z_k}^{-1}\cdots f_{Z_{j+1}}^{-1}\s_{Z_j}(Z_{j-1})$. These
divisors intersect transversally, and we define
$$D_J=\bigcap_{j\in J}D_j$$
for $J\subset[1,k]$. For such a subset $J$ of $[1,k]$ we can consider
the subword $\u_J=(\a_{i_j})_{j\not\in J}$ and there is a natural
isomorphism $\cFt_{\u_J}\simeq D_J$. We will therefore consider the
Bott--Samelson varieties $\cFt_{\u'}$ for any subword $\u'$ of $\u$ as
subvarieties of the Bott--Samelson variety $\cFt_\u$. We shall denote the
union of the divisors $D_j$ for $j\in[1,k]$ by $\partial\cFt_\u$.

Recall that by the construction, there is a map
$p_{\cFt_\u}:\cFt_\u\to\cF$. If $\u$ is a reduced expression for an
element $u$ of the Weyl group $W$, the natural map $p_{\cFt_\u}$ is
birational onto $\cF_u$ yielding a resolution of the Schubert
variety $\cF_u$.

\begin{rema}\label{rem-expression} 
The choice of a Bott--Samelson resolution $\cFt_\u$ for $\cF_u$ depends
on the choice of a reduced expression for $u$. Recall from
\cite{demazure} that since $\cF_w$ is a Schubert subvariety of
$\cF_v$, we may choose a reduced expression $\v=({\a_{i_1}}\cdots
\a_{i_k}\cdots \a_{i_l})$ for $v$ such that $\w=(\a_{i_k}\cdots
\a_{i_l})$ is a reduced expression for $w$. In particular in the
diagram
$$\xymatrix{\cFt_\w\ar@{^{(}->}[r]\ar[d]&\cFt_\v\ar[d]\\
\cF_w\ar@{^{(}->}[r]&\cF_v}$$
the vertical maps are birational and thus simultaneous resolutions of
singularities. We choose such a reduced expression $\v$ for $v$ to
construct $\cFt_\v$ and thus $\cFt_\w$.
\end{rema}

\subsection{Resolutions of $X$ and $Y$}
\label{resol}

\subsubsection{Two groups}

In this subsection, we will need to distinguish between the type $A$ and type $D$ cases.

In type $A$, set $G=\SL(\im N)$ and $G'=\SL(\K^n)$. We embed $G$ in $G'$ as
follows. Choose a complement $E_1$ for $\im N$ in $\ker N$ and
a complement $E_2$ for $\ker N$ in $\K^n$. We consider the
subgroup $G_0$ of $G'$ defined by:
$$G_0=\left\{f\in G'\ \left|
  \begin{array}{l}
    \textrm{$f$ stabilizes $\im N$, $\ker N$, and $E_i$} \ {\rm
  for}\ i\in\{1,2\},\ \\
\det(f\vert_{\im N})=\det(f\vert_{E_2})=1, \ {\rm
  and}\ f\vert_{E_1}=\id_{E_1}
\end{array}\right.\right\}.$$
This group is isomorphic to the product $\SL(\im N)\times \SL(E_2)$. Furthermore,
observe that $E_2$ is identified with $\K^n/\ker N$ and with $\im N$
via $N$. The group $G_0$ is thus isomorphic to $G\times G$; let us
embed $G$ into $G_0$ diagonally. For
any Borel subgroup $B$ of $G$, we may find a Borel subgroup $B'$ of $G'$
such that $B\subset B'$ in this embedding. We may thus consider the
variety $Z_0=G/B$ as a subvariety of $G_0/B_0$ and also of
$G'/B'=\cF$. We thus have a map $p_{Z_0}:Z_0\to\cF$.

In type $D$, set $G=\Sp(\im N)$ (recall that we have a non degenerate
skew form $\a$ on $\im N$) and $G'=\SO(\K^{2n})$. We embed $G$ in $G'$ as
follows. Choose an orthogonal complement $E_1$ for $\im N$ in $\ker N$ and
an isotropic subspace $E_2$ in $\K^{n}$ mapping bijectively to
$\K^{2n}/\ker N$. We consider the
subgroup $G_0$ of $G'$ defined by:
$$G_0=\left\{f\in G'\ \left|
  \begin{array}{l}
    \textrm{$f$ stabilizes $\im N$, $\ker N$, and $E_i$} \ {\rm
  for}\ i\in\{1,2\},\ \\
\det(f\vert_{\im N})=\det(f\vert_{E_2})=1, \ {\rm
  and}\ f\vert_{E_1}=\id_{E_1}
\end{array}\right.\right\},$$
remark that here the conditions $f(\ker N)\subset \ker N$ and
$\det(f\vert_{E_2})=1$ are
redundant. This group is isomorphic to $\SL(\im N)$.
The group $G$ embeds into $G_0$. For any Borel subgroup $B$ of $G$,
we may find a Borel subgroup $B'$ of $G'$ such that $B\subset B'$ in
this embedding. We may thus consider the variety $Z_0=G/B$ as a
subvariety of $G_0/B_0$ and also of $G'/B'=\cF$. We thus have a map
$p_{Z_0}:Z_0\to\cF$.

\subsubsection{Resolutions}

We again deal with types $A$ and $D$ simultaneously.

Let us take a sequence of simple roots
$\u=(\a_{i_1},\cdots,\a_{i_k})$ and apply the same
construction as for the Bott--Samelson variety $\cFt_\u$, but starting with
$Z_0=G/B$. We get a variety $\Xt_\u$ together with a morphism
$p_{\Xt_\u}:\Xt_\u\to\cF$. This variety can also be seen as the
homogeneous fiber bundle $\Xt_\u=G\times^B\cFt_\u$ where the action of $B$ on
$\cFt_\u$ is induced by the inclusion $B\subset B'$ and the natural
action of $B'$ on $\cFt_\u$. For any subword $\u'$ of $\u$, the
variety $\Xt_{\u'}$ can again be realized as a complete intersection
in $\Xt_\u$. In particular we have the same description of divisors
$G\times^B D_j$ for $j\in[1,k]$ on $\Xt_\u$ as on $\cFt_\u$. We
shall denote the union of these
divisors by  $\partial\Xt_\u$. Finally, we have a natural map $p_{\Xt_\u}:\Xt_\u\to\cF$.

Using the reduced expression $\v$ of $v$ defined in Remark
\ref{rem-expression}, we obtain a variety $\Yt=\Xt_\v$ and a subvariety
$\Xt=\Xt_\w$ of $\Yt$. We have natural maps 
$p_{\Xt_\v}$ (resp. 
$p_{\Xt_\w}$) from $\Yt$ (resp. $\Xt$) to $\cF$. Since the
maps $p_{\cFt_\v}$ and $p_{\cFt_\w}$ are $B'$-equivariant and thus
$B$-equivariant, we get a diagram
$$\xymatrix{\Xt=G\times^B\cFt_\w\ \ar@{^{(}->}[r]\ar[d]&
\Yt=G\times^B\cFt_\v\ar[d]\\
\Xh=G\times^B\cF_w\ar@{^{(}->}[r]\ar[rd]&\Yh=G\times^B\cF_v\ar[d]\\
&\cF.}$$
where the morphisms $G\times^B\cF_v\to\cF$ and
$G\times^B\cF_w\to\cF$ are given by $(g,x)\mapsto g\cdot x$.
The maps $p_\Yt$ and $p_\Xt$ are the vertical compositions in
the above diagram. We also have the projection maps $q_\Yt:\Yt\to G/B$ and
$q_\Xt\colon\Xt\to G/B$.

We must again make a distinction between type $A$ and type $D$ cases. Let us define $\cF'$ to be $\cF$ in type $A$ and ${\cal L}_\a(\im N)\times\Or\cF$ in type $D$. There are natural maps $p_\Yt$ and $p_\Xt$ from $\Yt$ and $\Xt$ to $\cF'$. In type $A$ they are simply defined by $p_\Yt=p_{\Xt_\v}$ and $p_\Xt=p_{\Xt_\u}$. In type $D$, the variety $G/B$ is isomorphic to $\Sp\cF(\im N)$ and there is a projection ${\rm pr}:G/B\to {\cal L}_\a(\im N)$. We thus have maps $\Yt\to {\cal L}_\a(\im N)$ and $\Xt\to {\cal L}_\a(\im N)$ defined by the composition of $q_\Yt$ and $q_\Xt$ with ${\rm pr}$. The maps $p_\Yt:\Yt\to\cF'={\cal L}_\a(\im N)\times\cF$ and $p_\Xt:\Xt\to\cF'={\cal L}_\a(\im N)\times\cF$ are the products of these maps with $p_{\Xt_\v}$ and $p_{\Xt_\u}$.

\begin{prop}
\label{prop-birat}
(\i) The maps $q_\Yt$ and $q_\Xt$ are dominant and locally trivial
fibrations with fiber over $(F_k)_{k\in[0,r]}\in G/B$ isomorphic to
Bott--Samelson varieties $\cFt_\v$ and $\cFt_\w$, respectively.

(\i\i) The maps $p_\Yt\colon \Yt\to \cF$ and $p_\Xt\colon\Xt\to\cF$ are
birational and dominant onto $Y$ and $X$ respectively. Thus they are resolutions of
singularities for $Y$ and $X$.
\end{prop}

\begin{proo}
The first part is clear from the definition of $\Yt$ and $\Xt$. The
second part follows from the birationality of the Bott--Samelson resolutions
$\cFt_\v\to\cF_v$ and $\cFt_\w\to\cF_w$, the smoothness of $\Yt$ and $\Xt$, and the first
part.
\end{proo}

\begin{nota}
  For  a subword $\u$ of $\v$, we define $X_\u$ to be the subvariety of
  $Y$ obtained as the image of $\Xt_\u$ (seen as a subvariety of
  $\Yt=\Xt_\v$) under the map $p_\Yt$. With this notation $X=X_\w$.
The map $p_\Xt:\Xt\to X$ is the resolution $\pi$ in Theorem \ref{main2}.
\end{nota}

\subsection{Existence of a splitting}

We have the following

\begin{theo}
\label{theo-split}
(\i) There exists a $B'$-canonical splitting of the Bott-Samelson variety
$\cFt_\v$ compatibly splitting all Bott-Samelson subvarieties
$\cFt_\u$ of $\cFt_\v$ for each subword $\u$ of $\v$.

(\i\i) This splitting induces a $B$-canonical splitting of $\Yt$
compatibly splitting all the subvarieties $\Xt_\u$ for $\u$ a subword of
$\v$.

(\i\i\i) The latter splitting induces a splitting of $Y$ compatibly splitting
all the subvarieties $X_\u$, where $\u$ is a subword of $\v$.
\end{theo}

\begin{proo}
  \textit{(\i)} This is an application of
\cite[Proposition 4.1.17]{BK}.

\textit{(\i\i)} We first observe that the $B'$-canonical splitting in \textit{(\i)} is a
$B_0$-canonical splitting, where $B_0=B'\cap G_0$ and $G_0$ was
defined in section \ref{resol}. For this, use the following result
(see \cite[Lemma
4.1.6]{BK}): let $H$ be a connected and simply connected semisimple group,
let $H'$ be a Borel
subgroup in $H$, and let $H''$ be a maximal torus in $H$. Let $X$ be a $H'$-scheme and
let $\phi\in\Hom(F_*\co_X,\co_X)$, where $F$
is the Frobenius morphism. Let us denote by $e_\a^{(n)}$ the divided powers,
where $\a$ is a root of $H$. There exists a natural action $e_\a^{(n)}*\phi$
of $e_\a^{(n)}$ on $\phi$ (see \cite[Definition 4.1.4]{BK}).

\begin{lemm}
The element $\phi$ is a $H'$-canonical splitting if and only if $\phi$
is $H''$-invariant and $e_\a^{(n)}*\phi=0$ for all $n\geq p$ and $\a$ a simple
root.
\end{lemm}

In our situation, we easily check that the divided powers of $G_0$ are divided
powers for $G'$. In particular the splitting in \textit{(\i)} is a $B_0$-canonical
splitting and compatibly splits the varieties $\cF_\u$.

To prove that the $B_0$-canonical splitting induces a $B$-canonical
splitting, we use the results of X. He and J.F. Thomsen \cite{he-thomsen}. More precisely we use Theorem 4.1 in that paper with $\lambda=\rho_{G_0}$. The hypothesis of that theorem are satisfied by Section 8 in \cite{he-thomsen}. We therefore obtain that $\Yt=G\times^B\cFt_\v$ has a Frobenius spliting which compatibly splits the subvarieties $G\times^B\cFt_\u$ for all subwords $\u$ of $\v$.\footnote{Note that in type $A$, \emph{i.e.} for a diagonal embeding, the result of X. He and J.F. Thomsen above was already proved by the same authors in \cite[Theorem 7.2]{HT}.}

We thus have a $B$-canonical splitting on $\cFt_\v$ compatible
with all the divisors $D_j$ and therefore with all the subvarieties
$\cFt_\u$ (by \cite[Proposition 1.2.1]{BK} and the fact that
the varieties $\cFt_\u$ are intersections of such divisors. Applying
\cite[Theorem 4.1.17]{BK}, we get a $B$-canonical splitting on
$G\times^B\cFt_\v=\Yt$ compatible with all subvarieties
$G\times^B\cFt_\u=\Xt_\u$.

\textit{(\i\i\i)} This is a direct application of Lemma 1.1.8 in \cite{BK}
together with the fact that $Y$ is normal (in type $A$ it is a Schubert variety and in type $D$ see Remark \ref{rema-norm}).
\end{proo}

\subsection{$D$-splitting}

In this subsection, we prove that the previous splitting is a $D$-splitting
for an explicit divisor $D$. For this we first need to compute
the canonical divisor of the variety $\Yt$.

Let us first fix some notation. As we have seen, if
$\v=(\a_{i_1},\cdots,\a_{i_k},\cdots,\a_{i_l})$ and if we denote by
$\v[j]$ the
subword $(\a_{i_1},\cdots,\a_{i_j})$ for $j\in[1,l]$, then the variety
$\Yt$ can
be realized as a sequence of $\pu$-fibrations $\Yt=\Xt_\v\to
\Xt_{\v[l-1]}\to\cdots\to \Xt_{\v[1]}\to G/B$. For all $j\in[1,l]$,
there is a natural map $p_{\Xt_{\v[j]}}:\Xt_{\v[j]}\to\cF$ and if $\co_\cF(1)$ is
the line bundle on $\cF$ defined by the Pl{\"u}cker embedding, we define
$\cL_{\v[j]}=p^*_{\Xt_{\v[j]}}(\cO_\cF(1))$. We shall denote by
$\cL_\Yt$ and $\cL_\Xt$ the line bundles $\cL_\v$ and $\cL_\w$,
respectively.

The following lemma is an easy modification of a well-known
result on the canonical divisor of the Bott--Samelson resolution,
see for example \cite[Proposition 2.2.2]{BK} or
\cite[Proposition 8.1.2]{kumar}:

\begin{lemm}
\label{cano}
We have the equality $\omega_{\Yt}^{-1}=
\co_{\Yt}(\partial \Yt)\otimes\cL_\Yt$.
\end{lemm}

\begin{proo}
  We prove the following formula by induction over $j\in[0,l]$:
$$\omega_{\Xt_{\v[j]}}^{-1}=\co_{\Xt_{\v[j]}}(\partial
\Xt_{\v[j]})\otimes\cL_{\v[j]}.$$
For $j=0$, we have $\Xt_{\v[j]}=G/B$. The line bundle $\omega_{G/B}^{-1}$
is the line bundle $\cL_{\v[j]}$ which is twice the
ample line bundle defined by the Pl{\"u}cker embedding of
$G/B$, since $G/B$ is diagonally embedded into $\cF=G'/B'$. Let us denote
the fibration $\Xt_{\v[j+1]} \to\Xt_{v[j]}$ by $f$  and its section by $\sigma$. The
induction follows from the formula:
$$\omega_{\Xt_{\v[j+1]}}=f^*\omega_{\Xt_{\v[j]}}\otimes
\co_{\Xt_{\v[j+1]}}(-\Xt_{\v[j]})\otimes \cL_{\v[j+1]}^{-1}\otimes
f^*\sigma^*\cL_{\v[j+1]}$$
which is a direct application of \cite[Lemma A-18]{kumar} and the fact that for each $j$ the divisor $\cL_{\v[j+1]}$ has the relative degree 1 for the fibration $f$. To prove the induction step we also the the equality  $\sigma^*\cL_{\v[j+1]}=\cL_{\v[j]}$.
\end{proo}

\begin{theo}
  There exists a $\Dt$-splitting of $\Yt$ compatibly splitting
  the subvarieties $\Xt_\u$, where $\u$ is a subword of $\v$ and $\Dt$ is an
  effective divisor such that $\cO_{\Yt}(\Dt)={\cL}_\Yt^{\otimes p-1}$.
\end{theo}

\begin{proo}
  Recall that in Theorem \ref{theo-split} we constructed a splitting
  $\varphi$
of $\Yt$ compatibly splitting the subvarieties $\Xt_\u$ for  a
subword $\u$ of $\v$. In particular, it is compatible with each of the divisors $\Xt_{\v(j)}$ for
$j\in[1,l]$, where
$\v(j)=(\a_{i_1},\cdots,\widehat{\a_{i_j}},\cdots,\a_{i_l})$. By
\cite[Theorem 1.4.10]{BK}, the splitting $\varphi$ provides a
$(p-1)\Xt_{\v(j)}$-splitting for all $j\in[1,l]$.
We may thus write
$${\rm div}(\varphi)=(p-1)\sum_{j=1}^l\Xt_{\v(j)}+\Dt=\partial\Yt+\Dt$$
with $\cO_{\Yt}(\Dt)=\cL_\Yt^{\otimes p-1}$ (compare with Lemma
\ref{cano}). But again by \cite[Theorem 1.4.10]{BK}, the
splitting $\varphi$ is a ${\rm div}(\varphi)$-splitting.
Now using \cite[Remark 1.4.2 (ii)]{BK} we get that $\varphi$ is a
$\Dt$-splitting. Note that by \cite[Theorem 1.3.8]{BK}, the
multiplicity of a divisor in the subscheme of zeros of the splitting
(here in ${\rm div}(\varphi)$) is at most $p-1$ thus $\Dt$ has a support
disjoint from any of the divisors $\Xt_{\v(j)}$. The divisors
$\Xt_{\v(j)}$ are therefore $\Dt$-compatibly split.

Finally, the result follows from the fact that for any subword
$\u$ of $\v$ the  variety $\Xt_\u$ is the intersection of
certain divisors $\Xt_{\v(j)}$.
\end{proo}

\begin{coro}
  There exists a $D$-splitting of $Y$ compatibly splitting all the
subvarieties $X_\u$ for a subword $\u$ of $\v$ for an effective divisor  $D$ such that
$\cO_{Y}(D)=({\cO_\cF(1)\vert_{Y}})^{\otimes p-1}$.
\end{coro}

\begin{proo}
  This is a direct application of \cite[Lemma 1.4.5]{BK} to the map
$p_\Yt:\Yt\to Y$. Indeed, in the previous Theorem, the divisor $\Dt$ is
the pullback by $p_\Yt$ of a divisor $D$ with the above
property. We may apply Lemma 1.4.5 in \cite{BK} because $Y$ is normal and
the conclusion follows for the splitting of the varieties $X_\u$ because
these varieties are the images of the varieties $\Xt_\u$ under $p_\Yt$.
\end{proo}

\begin{rema}
  The divisor $D$ may not be ample on $Y$ but its restriction to $X$
  is ample as $X$ is a subvariety of $\cF$.
\end{rema}

\section{Normality}

In this section we prove the results stated in the introduction. The
proof will be similar to the proof of the same results for Schubert
varieties as given in the book \cite{BK}. In order to pass from
positive characteristic to characteristic zero, we shall use the
results in \cite[Section 1.6]{BK}. For this we need to realize the
Springer fiber over $\Z$. This can be easily done by choosing a
representative of the nilpotent element $N$ in the normal Jordan
form in its $\GL(\K^n)$-orbit.

To simplify notation also also because the variety $Y=X_\v$ and some of the varieties $X_\u$ for $\u$ a subword of $\v$ are not contained in $\cF$ in type $D$ (but in $\cF'$), we will only consider in this section the varieties $X_\u$ for $\u$ a subword of $\w$. These varieties are subvarieties of $X=X_\w$ and can therefore be considered as subvarieties of $\cF$. In particular, the restriction of the above divisor $D$ is ample on these varieties.

\subsection{Some preliminary results}

We prove the normality of all the varieties $X_{\w[j]}$ for
$j\in[0,l-k]$ by induction over $j$. For this we need a more precise
description of the geometry relating $\Xt_{\w[j]}$ and
$\Xt_{\w[j+1]}$. Recall the construction of the variety
$\Xt_{\w[j+1]}$ from $\Xt_{\w[j]}$ by the elementary construction
\ref{cons-elementaire} as the fiber product
$\Xt_{\w[j+1]}=\Xt_{\w[j]}\times_\cF\cF_{\a_{i_{j+1}}}$. For a
subvariety $Z$ in $\cF$ we denote by $Z^{{i_{j+1}}}$ its image
under the
projection $p_{\a_{i_{j+1}}}:\cF\to\cF(\a_{i_{j+1}})$ (see Subsection
\ref{section-BS}).
We have the equality
$X_{\w[j+1]}=X_{\w[j]}^{i_{j+1}}\times_{\cF(\a_{i_{j+1}})}\cF$, so we
obtain the following
commutative diagram
\begin{equation}
  \label{diagramme}
\xymatrix{\Xt_{\w[j+1]}=\Xt_{\w[j]}\times_\cF\cF_{\a_{i_{j+1}}}
  \ar[r]^{\ \ \ \
  \ \ \ \ \ \ \ \ \widetilde{a}}
\ar[d]^{b'}&\Xt_{\w[j]}\ar[d]^{p_{\Xt_{\w[j]}}}\\
X'_{\w[j+1]}:=X_{\w[j]}\times_\cF\cF_{\a_{i_{j+1}}}\ar[r]^{\ \ \ \
  \ \ \ \ \ \ \ \ a'}\ar[d]^{b}&
X_{\w[j]}\ar@{_{(}->}[ld]\ar[d]^p\\
X_{\w[j+1]}=X_{\w[j]}^{i_{j+1}}\times_{\cF(\a_{i_{j+1}})}\cF\ar[r]^{\ \ \ \
  \ \ \ \ \ \ \ \ a}
\ar[d]^p&X_{\w[j]}^{i_{j+1}}\ar@{_{(}->}[ld]\\
X_{\w[j+1]}^{i_{j+1}}.&}
\end{equation}

\begin{lemm} \label{lemm-birat} With the above notation,

(\i) the map $p_{\Xt_{\w[j]}}:\Xt_{\w[j]}\to X_{\w[j]}$ is birational for
  all $j\in[0,l]$;

(\i\i) the map $p:X_{\w[j]}\to X_{\w[j]}^{i_{j+1}}$ is birational for
  all $j\in[0,l-1]$;

(\i\i\i) we have the equality
$X_{\w[j]}^{i_{j+1}}=X_{\w[j+1]}^{i_{j+1}}$ for
  all $j\in[0,l-1]$.
\end{lemm}

\begin{proo}
\textit{(\i)} We prove this by descending induction on $j$. For $j=l-k$ the
  corresponding map is $p_\Xt:\Xt\to X$, which is birational by
  Proposition \ref{prop-birat}.

Assume that $p_{\Xt_{\w[j+1]}}$ is birational. This map is the
  composition of the top two left vertical arrows $b$ and $b'$ in the
  previous diagram. In particular these two maps $b$ and $b'$ are also
  birational. But the topmost right vertical arrow in the above diagram
gives $b'$ by fiber product. This map is $p_{\Xt_{\w[j]}}$ and has to
  be birational.

\textit{(\i\i)} By what we just proved the map $b$ is birational. But it is a
fiber product of the map $p:X_{\w[j]}\to X_{\w[j]}^{i_{j+1}}$, which
has to be birational.

\textit{(\i\i\i)} Recall that we have two maps $\varphi_{\a_{i_{j+1}}}$ and
$\psi_{\a_{i_{j+1}}}$ from $\cF_{\a_{i_{j+1}}}$ to $\cF$ obtained by forgetting
one of the two subspaces corresponding to points in $\G(\a_{i_{j+1}})$
in $\cF_{\a_{i_{j+1}}}$. The
map $a'$ corresponds to $\varphi_{\a_{i_{j+1}}}$ while $b$ corresponds to
$\psi_{\a_{i_{j+1}}}$. The composition of the two forgetful maps
$\varphi_{\a_{i_{j+1}}}$ and $\psi_{\a_{i_{j+1}}}$ yield a map
$\cF_{\a_{i_{j+1}}}\to\cF^{i_{j+1}}$. The maps $p\circ a'$ and $p\circ b$
correspond by fiber product to the maps
$\psi_{\a_{i_{j+1}}}\circ\varphi_{\a_{i_{j+1}}}$
and $\psi_{\a_{i_{j+1}}}\circ\varphi_{\a_{i_{j+1}}}$, respectively. In
particular these two maps are equal and the result follows.
\end{proo}

\subsection{Proof of Theorem \ref{main1}}

We prove by ascending induction over $j$ that $X_{\w[j]}$ is normal. For $j=0$, we
have $X_{\w[0]}\simeq G/B$, which is normal. Let $j>0$ and assume that
$X_{\w[j]}$ is normal. The map
$a:X_{\w[j+1]}\to X_{\w[j]}^{i_{j+1}}$ is a $\pu$-fibration. Thus
to prove the normality of $X_{\w[j+1]}$ we only need to prove the
normality of $X_{\w[j]}^{i_{j+1}}$. But the map $p:X_{\w[j]}\to
X_{\w[j]}^{i_{j+1}}$ is birational and surjective (Lemma
\ref{lemm-birat} \textit{(\i\i)}), and $X_{\w[j]}$ is normal by the induction
hypothesis, so we only need to prove the equality
$$p_*\cO_{X_{\w[j]}}=\cO_{X_{\w[j]}^{i_{j+1}}}.$$
This will be done using the following lemma (see \cite[Lemma
3.3.3]{BK}):

\begin{lemm}
  Let $f\colon X\to Y$ be a surjective morphism between projective schemes,
  and let $\cL$ be an ample line bundle on $Y$. Assume that the map
  $H^0(Y,\cL^\nu)\to H^0(X,f^*\cL^\nu)$ is an isomorphism for $\nu$ large enough. Then $f_*\cO_X=\cO_Y$.
\end{lemm}

Consider $\cL$ ample on $X_{\w[j]}^{i_{j+1}}$ and the following
commutative diagram
$$\xymatrix{H^0(X_{\w[j+1]}^{i_{j+1}},\cL)\ar[r]\ar@{=}[d] &
  H^0(X_{\w[j+1]},p^*\cL)\ar[d]\\
H^0(X_{\w[j]}^{i_{j+1}},\cL)\ar[r]& H^0(X_{\w[j]},p^*\cL).\\
}$$
But since $X_{\w[j+1]}$ is $D$-split compatibly with $X_{\w[j]}$ and $D$
is ample, we get by \cite[Theorem 1.4.8]{BK} that the right
vertical map is surjective (the invertible sheaf $p^*\cL$ is semi-ample as the pull-back of an ample line bundle). Moreover, the map $p:X_{\w[j+1]}\to
X_{\w[j+1]}^{i_{j+1}}$ is a $\pu$-fibration, so the top horizontal map
is also surjective. We obtain that the lower horizontal map is
surjective. It is injective since the map $p:X_{\w[j]}\to
X_{\w[j]}^{i_{j+1}}$ is surjective. We may thus apply the previous
lemma and deduce the normality.

\begin{rema}
This proof works for $\K$ of positive characteristic, but relies only on
vanishing of cohomology and surjectivity of restrictions on cohomology results
which pass, by semi-continuity, to characteristic zero. The same proof
therefore works for $\chara(\K)=0$.
\end{rema}

\subsection{Proof of Theorem \ref{main2}}

Recall the definition of a rational morphism and of rational
singularities for $\chara(\K)=0$:

\begin{defi}
(\i) A morphism  of schemes $f\colon X\to Y$ is called \emph{rational} if
  $f_*\cO_X=\cO_Y$ and all its higher direct images vanish:
  $R^if_*\cO_X=0$ for $i>0$.

(\i\i) Assume $\chara(\K)=0$. A normal variety $X$ has \emph{rational
singularities} if there exists a rational birational proper morphism
$\pi:\Xt\to X$ with $\Xt$ smooth.
\end{defi}

We first prove the following

\begin{lemm}
\label{lemm-rat}
  For all $j\in[0,l-k]$, the map $p_{\Xt_{\w[j]}}:\Xt_{\w[j]}\to
    X_{\w[j]}$ is a rational morphism.
\end{lemm}

\begin{proo}
  We prove this lemma by induction over $j$. For $j=0$, we
have $p_{\Xt_{\w[0]}}:\Xt_{\w[0]}\to X_{\w[0]}$ is an isomorphism. Let
$j>0$ and assume that $p_{X_{\w[j]}}$ is rational. Then, since $b'$ is
obtained by fiber product from $p_{\Xt_{\w[j]}}$, we see that $b'$ is
rational. So, to prove the rationality of $p_{\Xt_{\w[j+1]}}$ we only
need to prove the rationality of $b$. Since $b$ is obtained by fiber
product from $p$, we only need to prove the rationality of
$p:X_{\w[j]}\to X_{\w[j]}^{i_{j+1}}$. But $p$ is birational
and $X_{\w[j]}^{i_{j+1}}$ is normal (see the proof of the normality of
$X$), therefore by the Zariski Main Theorem we obtain the equality
$p_*\cO_{X_{\w[j]}}=\cO_{X_{\w[j]}^{i_{j+1}}}$. We need to prove the
vanishing of the higher direct images. For this we embed $X_{\w[j]}$ in
$\cF$ and $X_{\w[j]}^{i_{j+1}}$ in
$\cF(\a_{i_{j+1}})$. We have a commutative diagram
$$\xymatrix{X_{\w[j]}\ar@{^{(}->}[r]\ar[d]^p&\cF\ar[d]^{p_{\a_{i_{j+1}}}}\\
X_{\w[j]}^{i_{j+1}}\ar@{^{(}->}[r]&\cF(\a_{i_{j+1}})}$$
and in particular the fibers of both morphisms are at most one-dimensional, thus $R^ip_*\cO_{X_{\w[j]}}=0$ for $i\geq 2$. Note that we also have the vanishing $R^2(p_{\a_{i_{j+1}}})_*\cE$ for any coherent sheaf $\cE$ on $\cF$. Thus the
surjection $\cO_\cF\to\cO_{X_{\w[j]}}$ induces a surjection
$R^1({p_{\a_{i_{j+1}}}})_*\cO_\cF\to R^1p_*\cO_{X_{\w[j]}}$. But since the second
vertical map is a $\pu$-fibration, we have $R^1({p_{\a_{i_{j+1}}}})_*\cO_\cF=0$ and the
result follows.
\end{proo}

Let us now prove the vanishing
${R^ip_{\Xt_{\w[j]}}}_*\omega_{\Xt_{\w[j]}}=0$ for $i>0$. For
this we use the following direct application of Theorem~1.2.12 from \cite{BK}:

\begin{lemm}
Let $f:X\to Y$ is a proper birational morphism with $X$ smooth. Assume
that $\varphi$ is a splitting for $X$ compatibly splitting a divisor
$Z$ such that the exceptional locus of $f$ is set-theoretically
contained in $Z$. Then we have $R^if_*\co_X(-Z)=0$ for $i>0$.
\end{lemm}

We want to apply this lemma to the map $p:\Xt_{\w[j]}\to X_{\w[j]}$,
the splitting constructed above, and the divisor $Z=\partial
\Xt_{\w[j]}$. For this we only need to check that the exceptional
locus of $p_{\Xt_{\w[j]}}$ is contained in $\partial \Xt_{\w[j]}$.
But the map $\Xt_{\w[j]}\to X_{\w[j]}$ decomposes as follows:
$$\Xt_{\w[j]}=G\times^B\cFt_{\w[j]}\to G\times^B\cF_{\w[j]}\to
X_{\w[j]}$$
and this map is $G$-equivariant. Furthermore, the complement to
$\partial\cFt_{\w[j]}$ in $\cFt_{\w[j]}$ is a $B$-equivariant dense
open subset thus the complement of $\partial \Xt_{\w[j]}$ in
$\Xt_{\w[j]}$ is a $G$-equivariant dense open subset. The map is
therefore an isomorphism on this open subset, and the exceptional locus
is contained in $\partial \Xt_{\w[j]}$. By the previous lemma, we get
the vanishing
$$R^i{p_{\Xt_{\w[j]}}}_*\co_{\Xt_{\w[j]}}(-\partial \Xt_{\w[j]})=0
\textrm{ for } i>0.$$
But from Lemma \ref{cano}, we have
$\omega_{\Xt_{\w[j]}}=\co_{\Xt_{\w[j]}}(-\partial \Xt_{\w[j]})\otimes
p_{\Xt_{\w[j]}}^*(\co_\cF(1)\vert_{X_{\w[j]}})$. Thus by projection
formula, we get:
$$R^i{p_{\Xt_{\w[j]}}}_*\omega_{\Xt_{\w[j]}} =
R^i{p_{\Xt_{\w[j]}}}_*\co_{\Xt_{\w[j]}}(-\partial \Xt_{\w[j]})\otimes
(\co_\cF(1)\vert_{X_{\w[j]}}) =0 \textrm{ for } i>0.$$
This completes the proof of Theorem~\ref{main2}. Corollary~\ref{main3}
follows from general results on rational resolutions, see
\cite[Lemma~3.4.2]{BK}, and Corollary~\ref{main4} follows from the
definition of rational singularities and Lemma~\ref{lemm-rat}.

\begin{rema}
The proof of Lemma~\ref{lemm-rat} works for any characteristic (once
the normality is proved). For $\chara(\K)=0$ we do not need to prove
the above vanishing $R^i{p_{\Xt_{\w[j]}}}_*\omega_{\Xt_{\w[j]}} =0$
for $i>0$. This result follows automatically from
Grauert--Riemenschneider Theorem \cite{GR}.
\end{rema}

\noindent
Nicolas {\sc Perrin}, \\
{\it Hausdorff Center for Mathematics,}
Universit{\"a}t Bonn, Villa Maria, Endenicher
Allee 62,
53115 Bonn, Germany, and \\
{\it Institut de Math{\'e}matiques de Jussieu,}
Universit{\'e} Pierre et Marie Curie, Case 247, 4 place
Jussieu, 75252 Paris Cedex 05, France.

\noindent {\it email}: \texttt{nicolas.perrin@hcm.uni-bonn.de}

\medskip\noindent
Evgeny {\sc Smirnov}, \\
{\it Department of Mathematics, Higher School of Economics,}
Myasnitskaya ul., 20, 101000 Moscow, Russia, and\\
{\it Laboratoire J.-V. Poncelet,}
Bolshoi Vlassievskii per., 11, 119002 Moscow, Russia.

\noindent {\it email}: \texttt{smirnoff@mccme.ru}

\end{document}